# Ruled Surfaces according to Rotation Minimizing Frame

Fatma Güler and Emin Kasap
Department of Mathematics, Arts and Science Faculty,
Ondokuz Mayis University, 55139 Samsun, Turkey.

**Abstract**

In this paper, we investigate the ruled surfaces generated by a straight line according to rotation minimizing frame (RMF). Using this frame of a straight line, we obtained the necessary and sufficient conditions when the ruled surface is developable. Also, we give some new results and theorems related to be the asymptotic curve, the geodesic curve and the line of curvature of the base curve on the ruled surface.

## 1. Introduction

A ruled surface is said to be "ruled" if it is generated by moving a straight line constinuously in Euclidean space $E^3$. Ruled surfaces are one of the simplest objects in geometric modeling.

One important fact about ruled surfaces is that they can be generated by straight lines. They can used in civil engineering. A special class of ruled surfaces, the developable surfaces, is characterized by a constant surface normal along each ruling. A surface is developable if and only if it is the envelope of a one-parameter family of planes. The developable surfaces are useful in geometric design, manufacturing systems and surface analysis. Also, they are used in computer aided design and manufacture.

The differential geometry of space curves is a classical subject which usually relates geometrical intuition with analysis and topology. Examination of special curves has an increasing interest in the theory of curves. One of the most significant surface curves is an asymptotic curve. The asymptotic curves pass through the directions of vanishing normal curvature on a surface. They are encountered in astronomy, astro-phssics and architectural CAD. Line of curvature is a useful tool in surface analysis for exhibiting variations of the principal direction. It can used in geometric design, shape recognition, polygonization of surfaces and surface rendering, [1-6]. Morever, line of curvature is closely related with rotation minimizing frames. Another of the most significant surface curves is a geodesic curve. A geodesic on a surface is the shortest path joining any two points of that surface. Geodesics are curves in surfaces that play a role anologous that of straight lines in the plane,[2].

The Frenet frame on the curve is a coordinate frame attached to the curve that helps describe the geometry of the curve. The useful property for curves parameterized by arc length is the Serret Frenet formulas. These formulas show the derivatives with respect to arc length of the Frenet frame as a function of the current Frenet frame, curvature and torsion. Another usable the frame of curve is rotation minimizing frame (RMF). They are characterized by the fact, that the normal plane of the curve rotates as little as possible around the tangent. They are useful in animation, motion planning swept surface constructions, and related applications where the Frenet frame may prove unsuitable,[8,9]. Also, rotation minimizing frames of space curves are used for sweep surface modeling in computer aided design.

Wang et al.[8] studied a novel simple and efficient method for accurate and stable computation of an RMF for any $C^1$ regular curve in 3D. This method called the double reflection method. Also, They analyzed the variational principles in design moving frames with boundary conditions, based on the RMF. In [4], the authors are proposed a new method to construct a developable surface possessing a given curve as the line of curvature of it. Also, they showed that the necessary and sufficient conditions when the resulting developable surface is a cylinder, cone or tangent surface. In[3], the authors analyzed the problem of rotional rotation minimizing frames (RMF). They addresses the rotation-minimizing osculating frame (RMOF) to construct ruled surfaces interpolating a space curve with tangent planes matching the osculating planes of that curve. In this paper, we obtained that the necessary and sufficient conditions the ruled surfaces generated by a straight line according to rotation minimizing frame (RMF) is developable. Also, using this frame of a straight line, we give some new results and theorems related to be asymptotic curve, geodesic curve and line of curvature of the base curve on the ruled surface.

## 2. Preliminaries

Let $r = r(s) : I \to E^3$ be a curve parametrized by the arc-length parameter $s$ with $r'(s) \neq 0$, where $r'(s) = \dfrac{dr(s)}{ds}$.

Useful property for curves parametrized by the arc-length is the Serret-frenet formulas. These formulas show the derivatives with respect to arc length of the Frenet frame as a function of the current Frenet frame, the curvature and the torsion. Then the Serret-frenet formulas of $r(s)$ are

$$
\begin{aligned}
T'(s) &= k(s) N(s) \\
N'(s) &= -k(s) T(s) + \tau(s) B(s) \\
B'(s) &= -\tau(s) N(s).
\end{aligned}
\qquad (2.1)
$$

where $T(s), N(s), B(s)$ are tangent, principal normal, and binormal vectors of the curve $r(s)$, the curvature $\kappa(s)$ and the the torsion $\tau(s)$ of the curve $r(s)$ at $s$.

Another usable the frame of curve is rotation minimizing frame. They are useful in animation, motion planning, swept surface constructions and related applications where the frenet frame may prove unsuitable. These frame is $\{T(s), U(s), V(s)\}$. The angular velocity vector $W$ for the frame $\{T(s), U(s), V(s)\}$ is defined by

$$W(s) = \beta(s) U(s) + \gamma(s) V(s).$$

That is $U$, $V$ are the rotation of $N$ and $B$ of the curve $r(s)$ in the normal plane. Then,

$$
\begin{bmatrix} U \\ V \end{bmatrix} = \begin{bmatrix} \cos\theta & \sin\theta \\ -\sin\theta & \cos\theta \end{bmatrix} \begin{bmatrix} N \\ B \end{bmatrix}
\qquad (2.2)
$$

where $\theta = \theta(s)$ represent the angle between the vectors $N$ and $U$.(see Fig.1), [4].

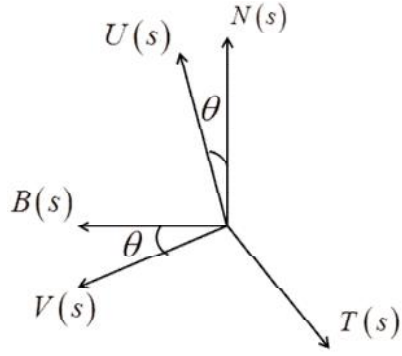

**Fig 1.** The Frenet frame $\{T(s), N(s), B(s)\}$ and the vectors $U(s)$, $V(s)$.

From Eqn. (2.2) differentiating the first-order derivatives of rotation minimizing frame may be expressed in the form as, [7]

$$U' = -\frac{\langle r'', U\rangle}{|r'|^2} r' , \quad V' = -\frac{\langle r'', V\rangle}{|r'|^2} r', \quad \theta' = -|r'|\tau . \tag{2.3}$$

The trace of an $\vec{X}$ oriented line along a space curve $r(s)$ is generally a ruled surface. A parametric equation of this ruled surface generated by $\vec{X}$ oriented line is given by

$$\phi(s,v) = r(s) + v\vec{d}(s), \quad s,v \in I \subset IR \tag{2.4}$$

where $\vec{d}$ is the unit direction vector of $\vec{X}$ oriented line.

**Definition 1.** The distribution parameter (or drall) of the $\phi(s,v)$ ruled surface is given as, by [5].

$$P_{\vec{d}} = \frac{\det(r'(s), \overrightarrow{d(s)}, \vec{d}'(s))}{\|\vec{d}'(s)\|^2}. \tag{2.5}$$

**Lemma 2.** The ruled surface Eqn. (2.4) is developable if and only if, by [5].

$$\det(r'(s), \overrightarrow{d(s)}, \vec{d}'(s)) = 0.$$

**Theorem 3.** A necessary and sufficient condition that a curve on a surface be a line of curvature is that the surface normals along the curve form a developable surface, [10].

**Definition 4.** For a curve $(r)$ lying on a surface, the following are well-known:

**(1)** $(r)$ curve is an asymptotic line of surface if and only if normal curvature $k_n$ vanishes.

**(2)** $(r)$ curve is an geodesic line of surface if and only if geodesic curvature $k_g$ vanishes, [5].

## 3. The ruled surfaces generated by a directed line according to rotation minimizing frame

Given a spatial curve $r: s \rightarrow r(s)$, where s is the arc length parameter. Let $\{T,N,B\}$ and $\{T,U,V\}$ be the Frenet frame field and rotation minimizing frame (RMF), respectively. We consider that a line $\vec{X}$ in $E^3$ such that it is firmly connected to rotating minimizing frame of $r(s)$ is represented, uniquely with respect to this frame in the form

$$X(s) = x_1(s)T(s) + x_2(s)U(s) + x_3(s)V(s) \tag{3.1}$$

where the components $x_i(s)$, $(i=1,2,3)$ are scalar functions of the arc length parameter of the base curve $r(s)$.

The trace of an $\vec{X}$ oriented line along a space curve $r(s)$ is a ruled surface. Such a surface has a parametrization in the ruled form as follows

$$\phi(s,v) = r(s) + v\vec{X}(s) \tag{3.2}$$

Differentiating Eqn. (3.1) with respect to $s$, we get Eqn. (2.3), then

$$\vec{X'}(s) = (x_1' - \kappa x_2 \cos\theta + \kappa x_3 \sin\theta)T + (\kappa x_1 \cos\theta + x_2')U + (x_3' - \kappa x_1 \sin\theta)V \tag{3.3}$$

From Eqn. (2.5) the distrubition parameter of these surface is,

$$P_X = \frac{(x_2 x_3' - x_3 x_2') - \kappa x_1(x_2 \sin\theta + x_3 \cos\theta)}{\sqrt{(x_1' - \kappa x_2 \cos\theta + \kappa x_3 \sin\theta)^2 + (\kappa x_1 \cos\theta + x_2')^2 + (x_3' - \kappa x_1 \sin\theta)}}$$

The unit normal vector to ruled surface $\phi(s,v)$ is given by

$$\overline{N} = \frac{\phi_s \times \phi_v}{\|\phi_s \times \phi_v\|} \tag{3.4}$$

Thus, from Eqn. (3.4) the unit normal vector to ruled surface $\phi(s,v)$ at the point $(s,0)$ is

$$\overline{N}(s,0) = \frac{x_2 V - x_3 U}{\sqrt{x_2^2 + x_3^2}} \tag{3.5}$$

The geodesic curvature of the base curve is

$$k_g = \langle \overline{N} \times T, T' \rangle = \frac{\kappa(x_2 \cos\theta - x_3 \sin\theta)}{\sqrt{x_2^2 + x_3^2}} \qquad (3.6)$$

and the normal curvature of the base curve is,

$$k_n = \langle \alpha'', \overline{N} \rangle = \frac{-\kappa(x_3 \cos\theta + x_2 \sin\theta)}{\sqrt{x_2^2 + x_3^2}}. \qquad (3.7)$$

And also the geodesic torsion is

$$\tau_g = \langle \overline{N} \times \overline{N}', T' \rangle = -\frac{\kappa^2}{x_2^2 + x_3^2}\left(\frac{1}{2}\sin 2\theta\left(x_3^2 - x_2^2\right) - \cos 2\theta\, x_2 x_3\right). \qquad (3.8)$$

Therefore, the following theorem may be given

**Theorem 1.** The base curve $r(s)$ is geodesic curve on the ruled surface $\phi(s,v)$ if and only if $\tan\theta = \dfrac{x_2}{x_3}$.

**Theorem 2.** The base curve $r(s)$ is asymptotic curve on the ruled surface $\phi(s,v)$ if and only if $\tan\theta = -\dfrac{x_3}{x_2}$.

**Theorem 3.** The base curve $r(s)$ is the line of curvature on the ruled surface $\phi(s,v)$ if and only if $\tan 2\theta = \dfrac{2 x_2 x_3}{x_3^2 - x_2^2}$.

**Corollary 5.** If the base curve $r(s)$ is geodesic curve and asymptotic curve on the ruled surface $\phi(s,v)$, then $\vec{X} = T$.

**Special Cases**

1. The case $\overrightarrow{X}(s) = T(s)$, $\overrightarrow{X}(s) = U(s)$, $\overrightarrow{X}(s) = V(s)$.

**i)** If $\overrightarrow{X} = T$, then $x_2 = x_3 = 0$, $x_1 = 1$. Thus from (3) $P_T = 0$.
**ii)** If $\overrightarrow{X} = U$, then $x_1 = x_3 = 0$, $x_2 = 1$. Thus from (3) $P_U = \theta' + \tau$.
**iii)** If $\overrightarrow{X} = V$, then $x_1 = x_2 = 0$, $x_3 = 1$. Thus from (3) $P_V = \theta' + \tau$.

**Corollary 6.** If $\overrightarrow{X} = T$, then from Lemma 2. the ruled surface $\phi(s,v)$ is developable. If $\overrightarrow{X} = U$ or $\overrightarrow{X} = V$, then the distibution parameters of the ruled surfaces $\phi(s,v)$ are equal.

2. The Case $\overrightarrow{X}(s) \in S_p\{T(s), U(s)\}$.

This case, we can write

$$\overrightarrow{X}(s) = x_1(s)T(s) + x_2(s)U(s). \tag{3.9}$$

Differentiating Eqn. (3.9) with respect to $s$, and Eqn. (2.3), then we get

$$\overrightarrow{X}' = (x_1' - \kappa x_2 \cos\theta)T + (x_2' + \kappa x_1 \cos\theta)U + (-\kappa x_1 \sin\theta)V. \tag{3.10}$$

From Eqn. (2.5) the distrubition parameter of the ruled surface $\phi(s,v)$ is

$$P_X = \frac{-\kappa x_1 x_2 \sin\theta}{(x_1' - \kappa x_2 \cos\theta)^2 + (x_2' + \kappa x_1 \cos\theta)^2 + (\kappa x_1 \sin\theta)^2}. \tag{3.11}$$

Thus, we may give the following corollary.

**Corollary 7.** The ruled surface $\phi(s,v)$ is developable if and only if $\theta = \pi + 2k\pi$, $k \in IR$.

3. The Case $\overrightarrow{X}(s) \in S_p\{T(s), V(s)\}$.

This case, we can write,

$$\overrightarrow{X}(s) = x_1(s)T(s) + x_3(s)V(s). \tag{3.12}$$

Differentiating Eqn. (3.12) with respect to $s$, and Eqn. (2.3), then we get

$$\overrightarrow{X}' = (x_1' + \kappa x_3 \sin\theta)T + (\kappa x_1 \cos\theta)U + (x_3' - \kappa x_1 \sin\theta)V \tag{3.13}$$

From Eqn. (2.5) the distrubition parameter of the ruled surface $\phi(s,v)$ is

$$P_X = \frac{-\kappa x_1 x_3 \cos\theta}{(x_1' + \kappa x_3 \sin\theta)^2 + (\kappa x_1 \cos\theta)^2 + (x_3' - \kappa x_1 \sin\theta)^2}. \tag{3.14}$$

Thus, we may give the following corollary.

**Corollary 8.** The ruled surface $\phi(s,v)$ is developable if and only if $\theta = \frac{\pi}{2} + 2k\pi$, $k \in IR$.

### 4. The Case $\vec{X}(s) \in S_p\{U(s), V(s)\}$.

Therefore, we can write,

$$\vec{X}(s) = x_2(s)U(s) + x_3(s)V(s). \tag{3.15}$$

Differentiating Eqn. (3.15) with respect to $s$, and Eqn. (2.3), then we get

$$\vec{X}' = (\kappa x_3 \sin\theta - \kappa x_2 \cos\theta)T + x_2'U + x_3'V \tag{3.16}$$

From Eqn. (2.5) the distrubition parameter of the ruled surface $P(s,t)$ is

$$P_X = \frac{x_2 x_3' - x_3 x_2'}{(\kappa x_3 \sin\theta - \kappa x_2 \cos\theta)^2 + (x_2')^2 + (x_3')^2}. \tag{3.17}$$

Thus, we may give the following corollary.

**Corollary 9.** The ruled surface $\phi(s,v)$ is developable if and only if $x_2 x_3' = x_3 x_2'$.

**Example 1..** Let $r(s) = \left(\frac{3}{5}\cos(s), \frac{3}{5}\sin(s), \frac{4}{5}s\right)$ be a unit speed curve. Then, RMF is easy to show that

$$\begin{vmatrix} T(s) = \left(-\frac{3}{5}\sin(s), \frac{3}{5}\cos(s), \frac{4}{5}\right), \\ U(s) = \left(\frac{4}{5}\cos(s)\sin\theta(s) - \cos\theta(s)\cos(s), -\sin(s)\cos\theta(s) - \frac{4}{5}\sin\theta(s)\cos(s), \frac{3}{5}\sin\theta(s)\right), \\ V(s) = \left(\cos(s)\sin\theta(s) + \frac{4}{5}\sin(s)\cos\theta(s), \sin(s)\sin\theta(s) - \frac{4}{5}\cos(s)\cos\theta(s), \frac{3}{5}\cos\theta(s)\right). \end{vmatrix}$$

If we take
$x_1(s) = s^2$, $x_2(s) = s^2$, $x_3(s) = s$ and $\tan\theta(s) = s$ then the Eqn. (3.18) is satisfied. Thus, we obtain geodesic curve $r(s)$ as

$$\phi(s,v) = r(s) + v(x_1(s)T(s) + x_2(s)U(s) + x_3(s)V(s)), \qquad (3.18)$$

where $-5 \leq s \leq 5$, $-1 \leq v \leq 1$ (Fig. 1).

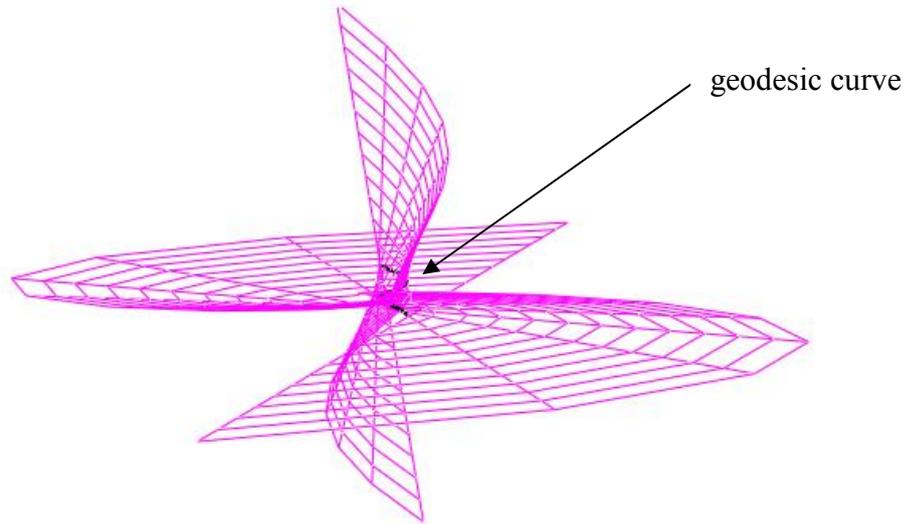

Fig.1. $r(s)$ is geodesic curve on the ruled surface

Let $r(s) = \left(\dfrac{3}{5}\cos(s), \dfrac{3}{5}\sin(s), \dfrac{4}{5}s\right)$ be a unit speed curve. Then RMF is easy to show that

$$\begin{cases} T(s) = \left(-\dfrac{3}{5}\sin(s), \dfrac{3}{5}\cos(s), \dfrac{4}{5}\right), \\ U(s) = \left(\dfrac{4}{5}\cos(s)\sin\theta(s) - \cos\theta(s)\cos(s), -\sin(s)\cos\theta(s) - \dfrac{4}{5}\sin\theta(s)\cos(s), \dfrac{3}{5}\sin\theta(s)\right), \\ V(s) = \left(\cos(s)\sin\theta(s) + \dfrac{4}{5}\sin(s)\cos\theta(s), \sin(s)\sin\theta(s) - \dfrac{4}{5}\cos(s)\cos\theta(s), \dfrac{3}{5}\cos\theta(s)\right). \end{cases}$$

If we take
$x_1(s) = s^2$, $x_2(s) = s$, $x_3(s) = -s^2$ and $\tan\theta(s) = s$ then the Eqn. (3.19) is satisfied.
Thus, we obtain asimptotik curve $r(s)$ as

$$\phi(s,v) = r(s) + v(x_1(s)T(s) + x_2(s)U(s) + x_3(s)V(s)), \qquad (3.19)$$

where $-5 \leq s \leq 5$, $-1 \leq v \leq 1$ (Fig. 2).

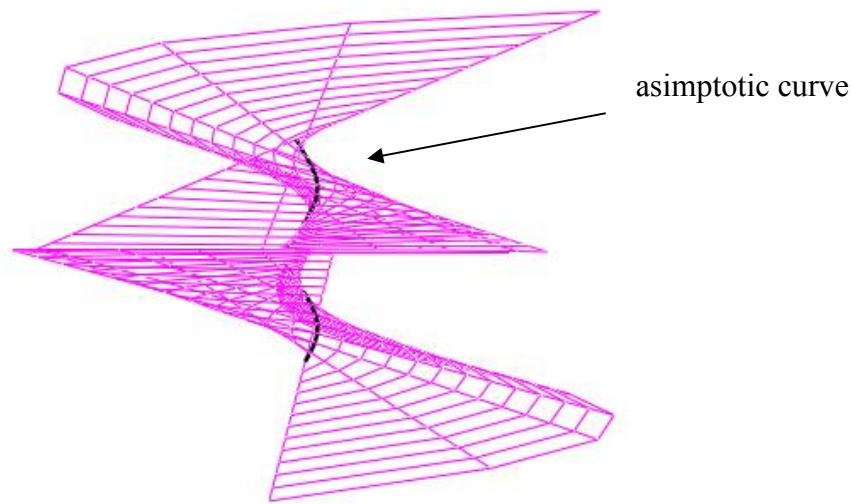

Fig.2. $r(s)$ is asimptotic curve on the ruled surface